\documentclass[12pt]{article}

\usepackage{amsmath} \usepackage{amsfonts} \usepackage{amssymb}

\newtheorem{defi}{Definition}
\newtheorem{theorem}{Theorem}
\newtheorem{cor}{Corollary}
\newtheorem{lemma}{Lemma}
\newtheorem{proposition}{Proposition}

\newcommand{\CC}{\mathbb{C}}
\newcommand{\RR}{\mathbb{R}}

\begin{document}
\title{ \bf Cohomology ring of symplectic quotients 
\\ by circle actions}
\author{Ramin Mohammadalikhani}

\date{}

\maketitle

\begin{abstract}

In this article we are concerned with how to compute the cohomology ring
of a symplectic quotient by a circle action using the information we have
about the cohomology of the original manifold and some data at the fixed
point set of the action.
Our method is based on the Tolman-Weitsman theorem which 
gives a characterization of the kernel of the Kirwan map.
First we compute a generating set for the kernel of the Kirwan map 
for the case of product of compact connected manifolds such that the
cohomology ring of each of them is generated by a degree two class.
We assume the fixed point set is isolated; however the circle action 
only needs to be ``formally Hamiltonian''. 
By identifying the kernel, we obtain the 
cohomology ring of the symplectic quotient.
Next we apply this result to some special cases and in particular 
to the case of products of two dimensional spheres. 
We show that the results of Kalkman
and Hausmann-Knutson are special cases of our result.

\end{abstract}

\section{ Introduction}

In this article we are concerned with the cohomology ring of
symplectic reductions. We would like to answer the following 
question:
When we consider a Hamiltonian action of a Lie group on a 
symplectic manifold, what would the quotient space topologically 
look like?
The interesting point is that in fact using only the information
about the moment map at the fixed point set of the
action one can at least theoretically answer this question.
The Tolman-Weitsman theorem \cite{TW1} has now 
enabled us to find the answer to our question with just 
the information mentioned .
Kalkman was the first who in \cite{Ka} calculated the cohomology
ring of the symplectic reduction of a projective space by a 
circle action using the localization formula.
However his work was not continued further.
The next attempt to understand the cohomology of these spaces 
was based on other means.

Hausmann and Knutson used Danilov's theorem to approach the 
problem.
Danilov's theorem specifies the cohomology rings of all 
toric manifolds.
In cases where the original manifold is a toric manifold, one can
show that its symplectic quotient is a toric manifold too.
One is then able to use Danilov's theorem to find the
cohomology ring of the symplectic reduction.
This is what Hausmann and Knutson did in
\cite{HK} to calculate the cohomology ring of the
symplectic quotient of a product of two-dimensional spheres.

We know very little when the group acting on the manifold is a general compact 
Lie group or even a torus.
In her Ph.D. thesis R. Goldin \cite{G} ( also see \cite{G2} )
answered this question for the action of a torus on a coadjoint
orbit of $ SU(n) $.

Later Tolman and Weitsman \cite{TW2} generalized the results of
Hausmann and Knutson to a compact connected symplectic manifold
but they had to assume that the action is semi-free and the
fixed point set is isolated.
They found the integer cohomology ring of
$ M_{red} = \mu^{-1}(0) / S^{1} $,
whenever $0$ is a regular value of the moment map.
The conditions of semi-free action and the fixed point set being
isolated enabled them to establish a correspondence between the
fixed point set of the circle action on $M$ and that of the 
product of two-dimensional spheres.
However to obtain those results Tolman and Weitsman did not use
their own theorem \cite{TW1}, which already opens the way to 
answer the problem in more general settings.
Our method in this article is based on this key theorem. 
This theorem reduces the problem of finding a generating set 
for the kernel of the Kirwan map 
$$ \kappa : \quad H^{*}_{T}(M) \quad \longrightarrow \quad
H^{*}(\mu^{-1}(0) / T)$$
to some specific algebraic calculations. We are then done 
with the task of finding the cohomology ring of the quotient 
space because
$$ H^{*}(\mu^{-1}(0) / T) \cong  H^{*}_{T}(M) / ker(\kappa), $$
due to Kirwan's surjectivity theorem (\cite{Ki1}, 5.4) which
asserts that $ \kappa $ is a surjective ring homomorphism.

We would like to state the Tolman-Weitsman theorem here for the 
case of circle actions on which the results of this article are 
based:

\begin{theorem}
\label{tw1}
\cite{TW1}

Let  $S^{1}$  act on a compact symplectic manifold M with moment 
map $ \mu : M \rightarrow \RR $.
Assume that $r$ is a regular value of $ \mu $. Let $ \mathcal{F}$
denote the set of fixed points of the action. Write
$ M_-(r) = \mu^{-1}(-\infty, r) \quad and \quad M_+ = 
\mu^{-1}(r,+\infty).$ 
Define
$K_{\pm}(r) = \{ \alpha \in H^{*}_{S^1}(M,\CC) : 
\alpha \mid_{\mathcal{F}\cap M_{\pm}(r)} = 0 \}  $
and
$ K(r) = K_+(r) \oplus K_-(r) .$
Then there is a short exact sequence:
$$0 \rightarrow K(r) \rightarrow   H^{*}_{S^1}(M,\CC) 
\stackrel{\kappa}
{\rightarrow} H^{*}(M_{red},\CC) \rightarrow 0, $$
where  $\kappa$  is the Kirwan map and
$ M_{red} = \mu^{-1}(r) / S^1 $.
\end{theorem}
When $r=0$ we write $ M_{\pm}$,$K_{\pm}$ and $K$ for 
 $ M_{\pm}(0)$,$K_{\pm}(0)$ and $K(0)$.

{\bf Remark:} By Remark 3.4 of \cite{TW1},
we do not need to assume that the action  is Hamiltonian.
The statement still holds if the action is more generally
{\it formally Hamiltonian}. This means there is a Morse-Bott
function  $ \mu: M \rightarrow  {\bf t}^* = Lie(S^1)^* \cong \RR $
(a formal moment map) 
such that the critical points of $ \mu $ correspond exactly to 
the fixed points of the action. Then as long as $M$ is compact 
and $0$ is a regular
value of $\mu$, the theorem is true for any formal moment map.
We still need to assume that $M$ is compact and $0$ is a regular 
value of the moment map.

Besides the Tolman-Weitsman theorem that we use in this article
the residue formula ( \cite{JK1}, \cite{JK3} ) is another 
powerful tool which may enable us to answer the question 
even in more general cases.

For the various definitions and properties of equivariant 
cohomology see for example \cite{Au} and \cite{BGV}.

\section{ Notation and Preliminaries}

First let us fix our notation.
Consider a compact connected manifold $M$ whose cohomology ring 
is generated by degree two classes
$ x_i \in H^2(M);  i = 1, 2, ..., m$.
Assume the manifold is equipped with a circle action with isolated
fixed points.
We label the fixed point set by
$ F_j ; j = 1, 2, ..., n$.
Suppose there are moment maps for the action denoted by
$ \mu_i : M \longrightarrow \RR $ such that 
$ i_t x_i = t d \mu_i $
for all $t \in \RR \cong {\mathbf t}^* =  Lie(T)^* $.
Here $T = S^1$. 
Consider the two-form $ x = \sum_{i=1}^m x_i $.
Corresponding to this two-form we also have the function
$ \mu : M \longrightarrow \RR $ defined by 
$ \mu = \sum_{i=1}^m \mu_i $ so that
$ i_t x =  t  d \mu $.

We impose the extra condition that $\mu$ does not vanish at any
of the fixed points.

Now consider the equivariant cohomology algebra $ H_T^* (M)$. 
As a vector space it can be written as
$$ \mathcal{R} = H_T^* (M)  \cong H^*_T(point) \otimes H^*(M), $$
where
$ H^*(M) \cong \CC [x_1,...,x_m ] / \mathcal{I}.$
Here $ \mathcal{I} $ is the set of relations in $H^*(M)$.
Also $ H^*_T(point) = \CC [t] $, the polynomial ring in
the variable $t$.

If $ \tilde{x}_i = x_i + t \mu_i $ are the equivariant
extensions of the corresponding $x_i$'s, then we see that
$ \tilde{x}_1, ..., \tilde{x}_m $ together with $t$
generate the equivariant cohomology $  H_T^* (M) $ as a ring.
We also consider the equivariant
extension $ \tilde{x}  = x + t \mu   $.

The values of the moment maps at the fixed points are of great
importance.
We denote them as follows:
$ \mu_i (F_j) = \theta_{ij} $ so that the restriction of
$\tilde{ x}_i $ to the j-th component of the fixed point set is
$ \theta_{ij} $:
$ \tilde{x}_i \mid_{F_j} = \theta_{ij}t $. Then
$ \mu (F_j) = \sum_{i=1}^m \theta_{ij} $ and
$ \tilde{x} \mid_{F_j} = \sum_{i=1}^m \theta_{ij}t $.

Now we would like to specify $K_+$ and $K_-$ in the ring
$ \mathcal{R}$.
According to the Tolman-Weitsman theorem,
$$ K_+ = \{ \alpha \in H_T^* (M) : \alpha \mid_{F_j} = 0
\textup{ for all } 
j \textup{ such that } \mu(F_j) > 0 \}. $$
Equivalently,
$$ K_+ = \{ \alpha \in \mathcal{R} : \alpha 
( \theta_{1j}t, ...,\theta_{mj}t)
= 0 \textup{ for all }j  \textup{ such that } 
\mu(F_j)= \sum_{i=1}^m \theta_{ij}  > 0 \}. $$
The ideal $ K_- $ is defined similarly with the difference that $ > $ 
is replaced with $ < $ in the definition of the set.
We can consider $ K_+$ and $K_-$  as the intersection of a finite
number of ideals as follows: Consider the multivariable polynomial ring
$$ \bar{\mathcal{R}} = 
\CC [t][ \tilde{ x}_1, ..., \tilde{ x}_m ], $$
in the variables $ \tilde{ x}_i $ with coefficients in $ \CC[t] $
(the polynomial ring in one variable $t$ with complex 
coefficients).
Thus $H^*_T (M) $ is the quotient of $  \bar{\mathcal{R}} $ by an 
ideal of relations.

For $ 1 \leq j \leq n $ define the ideals
$$ \mathcal{I}_j = \{  \alpha \in H_T^* (M) : \alpha 
\mid_{F_j} = 0 \} \cong \{ \alpha \in  \bar{\mathcal{R}} :
\alpha ( \tilde{ x}_1 = \theta_{1j}t, ...,\tilde{ x}_m = \theta_{mj}t) = 0 \} $$
in $  \bar{\mathcal{R}} $. Then $ \bar{K}_+$ is the intersection 
of those $ \mathcal{I}_j $'s
that correspond to the $ j $'s for which the value of the moment 
map
$\mu$ is positive:
$$  \bar{K}_+ = \bigcap_{ 1 \leq j \leq n : \mu(F_j) > 0 } 
\mathcal{I}_j \qquad
\textup{ (similarly }  \quad \bar{K}_- = 
\bigcap_{ 1 \leq j \leq n : \mu(F_j) < 0 } \mathcal{I}_j ). $$

In fact we know the generators of each $ \mathcal{I}_j $.
They are simply $ \tilde{x}_1 - \theta_{1j}t $, ...,
$ \tilde{x}_m - \theta_{mj}t $.
The problem of classifying the intersection ideal 
(say by specifying a generating set) is very hard and still open! 
We can solve this problem for a special case 
that is important to our geometric concerns. In the next section we will
explain this special case and will show that a generating set for the intersection
 ideal exists such that each of its elements is a product of proper linear forms. 

\section{ The main result and its proof}

We shall consider the special case when $ M $ is a product of 
compact connected symplectic manifolds $ M_i, \quad i = 1,2,...,m $,
i.e. $ M = M_1 \times M_2 \times ... \times M_m $
such that the cohomology ring of each $M_i$
is generated by a two-form $ x_i \in H^2(M_i) $, i.e. 
$  H^* (M_i) = < x_i > $.
Consider the extensions of these forms to $M$ by
$ x_i = 1 \otimes ... \otimes x_i \otimes ... \otimes 1
 \in H^*(M) = \otimes_{i=1}^m H^*( M_i). $
Each $ M_i$ is equipped with a Hamiltonian circle action with
isolated fixed points. Consider the diagonal action on $M$.
The fixed points are labeled by $m$-tuples
$  \mathbf{F}  = ( F_{1j_1}, F_{2j_2}, ..., F_{mj_m}) $ 
for all choices of $ 1 \leq j_i \leq n_i $, where $ n_i$ is 
the number of the fixed points of $M_i$ with distinct moment
 map value. 
Here $F_{ij}$ denotes the union of those fixed points of $M_i$ 
whose value under the moment map $ \mu_i$ is $\theta_{ij}$. 
Therefore
$ j = j' \Longleftrightarrow \theta_{ij} = \theta_{ij'}$
 for all $ i, j, j' $.
If the value of  $ \mu_i $ at $ F_{ij} $ is denoted
by $ \theta_{ij} $, then
$ \mu ( \mathbf{F}) = \mu ( F_{1j_1}, F_{2j_2}, ... , F_{mj_m} ) 
= \sum_{i=1}^m \theta_{ij_i} $.  
The restrictions of each $ \tilde{x}_i$ and $ \tilde{x}$ to the 
fixed point 
$ \mathbf{F}  = ( F_{1j_1}, F_{2j_2}, ... , F_{mj_m} ) $
are given by
$\tilde{x}_i \mid_{\mathbf{F}} = \theta_{ij_i}t $
and
$\tilde{x}\mid_{\mathbf{F}} = \sum_{i=1}^m \theta_{ij_i}t $.

As usual we are concerned about the kernel of the Kirwan map:
$ K = K_+ \oplus K_- $. 
The following proposition is of fundamental importance to us:

\begin{proposition}
\label{linearterms}
The ideal $  \bar{K}_+ $ has a generating set such that each generator 
is a product of linear forms of the form
$ \tilde{x}_i - \theta_{ij_i}t $. 
Moreover the linear terms that appear in each generator are 
mutually distinct.
  The same as for $  \bar{K}_+ $ is true for $  \bar{K}_- $.
\end{proposition}

Note that here $i$ indexes the manifolds $M_i$ and $j_i$ indexes 
the fixed point set of the $i$-th manifold $M_i$.

The proposition is an immediate consequence of the following lemma.

\begin{lemma}
\label{algebraiclemma2}
Consider the ring
$ \bar{\mathcal{R}} = \CC [t][ \tilde{x}_1 , ..., \tilde{x}_m ],$
and consider the following finite set  
$ \mathcal{F}' = \{ F = ( \theta_{1j_1}t, ..., \theta_{mj_m}t ) 
\in  \CC [t]^m :  1 \leq j_i \leq n_i  \textup{ such that }
\theta_{ij} > \theta_{ij'} \textup{ for } j <j' \}$
of points in $ \CC [t]^m $,
where the real numbers $ \theta_{ij}$ and positive integers $ n_i $
are given.\\
Define 
$ \bar{\mathcal{F}}_+ = \{ 
( \theta_{1j_1}t, ..., \theta_{mj_m}t ) \in
\mathcal{F}'  :  \sum_{i=1}^m \theta_{ij_i} > c \}, $
where $c$  is some fixed real number.
It is a subset of $ \mathcal{F}' $: Let 
$ \mathcal{I}_+ = \{ \alpha \in  \bar{\mathcal{R}} : 
\alpha (F) = 0
\textup{ for all } F \in  \bar{\mathcal{F}}_+ \}. $
Then the ideal $ \mathcal{I}_+ $ has a generating set consisting of
polynomials each of which is a product of terms of the form
$ \tilde{x}_i - \theta_{ij_i} t $
which we will refer to as linear terms from now on.
The linear terms in each generator are mutually distinct.
If we replace the condition 
$  \sum_{i=1}^m \theta_{ij_i} > c $ with
$  \sum_{i=1}^m \theta_{ij_i} < c $, 
the statement is still true.
\end{lemma}

To prove this, we need the following algebraic lemma:
\begin{lemma}
\label{algebraiclemma1}
If $P(x_1,...,x_n) $ is a polynomial and $P(a_1,...,a_n) = 0$,
then there are polynomials $ Q_1, ..., Q_n$ in $x_1,...,x_n$ such
that
$ P = (x_1 - a_1) Q_1 + ... +  (x_n - a_n) Q_n. $
\end{lemma}

{ \bf Proof of Lemma \ref{algebraiclemma1}:} 
Since $P(a_1,...,a_n) = 0$, the Euclidean Lemma
tells us that there are polynomials $ Q_1$ and $ R_1(x_2,...,x_n) $
such that $  P = (x_1 - a_1) Q_1 + R_1. $
Then  $R_1(a_2,...,a_n) = 0$. Thus, we can proceed by induction on $k$ for
$R_k(a_k,...,a_n) $, to finally obtain
$ P = (x_1 - a_1) Q_1 + ... +  (x_n - a_n) Q_n + R_n. $
where $R_n$ is just a number. Then from  $P(a_1,...,a_n) = 0$, we get
$R_n = 0$. This completes the proof of the lemma.
$\hspace{1cm} \stackrel{\infty}{\smile}$

{ \bf Proof of Lemma \ref{algebraiclemma2}:}
The proof is by induction on $m$.
To understand how the induction works, we initially discuss 
both cases $ m=1 $ and $ m=2 $, even though mathematically
we only need to check the case $ m = 1 $.

So assume $ m = 1 $. We show $ \mathcal{I}_+ $ is generated by 
one element, i.e.
$ \prod_{ j : \theta_j > c } ( \tilde{x} -\theta_j t ). $
To see this, notice that 
$ \alpha ( \tilde{x} ) \in \mathcal{I}_+ $
if and only if
$ \alpha (F) = 0 $
for every 
$ F \in  \bar{\mathcal{F}}_+ $.
This means
$ ( \tilde{x} -\theta_j t ) $ divides $ \alpha $ for each $j$
with $ \theta_j > c $.
So their product also divides $ \alpha $, which is what we
wanted to prove.

Now consider the case $ m= 2 $ so that  
$  \bar{\mathcal{R}} =  \CC [t][ \tilde{x}_1, \tilde{x}_2 ]. $
We arrange the points of $ \bar{\mathcal{F}}_+ $ in the following table:
$$ \hspace{-4cm} ( \theta_{11}t, \theta_{21}t ), \qquad 
( \theta_{12}t,\theta_{21}t ), \quad
...  \quad ,( \theta_{1l}t,\theta_{21}t )$$
$$ \hspace{-4cm}...\hspace{2cm}...\hspace{3cm}...$$
$$  \hspace{-3.8cm} ( \theta_{11}t,\theta_{2m_l}t ), \quad 
( \theta_{12}t,\theta_{2m_l}t ), \quad ... \quad ,
( \theta_{1l}t,\theta_{2m_l}t )$$
$$\hspace{-7.2cm}... \hspace{2cm}...$$
$$ \hspace{-7cm}( \theta_{11}t ,\theta_{2m_2}t ), \quad 
( \theta_{12}t,\theta_{2m_2}t )$$
$$\hspace{-9.4cm}...$$
$$ \hspace{-9.4cm} ( \theta_{11}t,\theta_{2m_1}t ).$$

Here the $k$-th row and $i$-th column is
$ ( \theta_{1i}t, \theta_{2k}t ) $ for $ k \leq m_i $.
The integer $ m_i $ is the largest integer such that
$ \theta_{1i} + \theta_{2m_i} > c $ and 
$ l $ is the largest integer such that
$ \theta_{1l} + \theta_{21} > c $. Notice that $ l \leq n_1 $.

Since
$ \theta_{11} > \theta_{12} > ... > \theta_{1n_1}$ and
$ \theta_{21} > \theta_{22} > ... > \theta_{2n_2}, $
then if $ \theta_{1(i+1)} + \theta_{2j} > c $, we also have
 $ \theta_{1i} + \theta_{2j} > c $. Therefore  
$ m_1 \geq m_2 \geq ... \geq m_l $
which is a crucial fact in our argument. 

Fix $ \alpha \in \mathcal{I}_+ $.
Then $ \alpha( \tilde{x}_1, \tilde{x}_2 ) $ vanishes at all
$ F \in  \bar{\mathcal{F}}_+$. 
By Lemma \ref{algebraiclemma1} applied to the first point of
the first column, we see that there are polynomials
$ p ( \tilde{x}_1, \tilde{x}_2 ) $ and $ q ( \tilde{x}_2 ) $
such that
$ \alpha( \tilde{x}_1 , \tilde{x}_2) =  ( \tilde{x}_1 - \theta_{11}t )
p ( \tilde{x}_1 , \tilde{x}_2 ) +  ( \tilde{x}_2 - \theta_{21}t )
q ( \tilde{x}_2 ) $.
Since $ \alpha $ vanishes at other points of the first column,
we see that
$ q ( \theta_{22}t) = ... = q ( \theta_{2m_1}t ) = 0 $
so that
$ ( \tilde{x}_2 - \theta_{22}t ) ( \tilde{x}_2 - \theta_{23}t )
 ...( \tilde{x}_2 - \theta_{2m_1}t ) $
has to divide  $ q ( \tilde{x}_2 ) $.
Therefore, there is a polynomial $ q' ( \tilde{x}_2 ) $
such that
$  q ( \tilde{x}_2 ) = ( \tilde{x}_2 - \theta_{22}t )
( \tilde{x}_2 - \theta_{23}t )
 ...( \tilde{x}_2 - \theta_{2m_1}t )  q' ( \tilde{x}_2 ) $.
Now by considering the vanishing of $\alpha$ at the first point of
the second column we find that there are polynomials
$ p_1 ( \tilde{x}_1, \tilde{x}_2 ) $ and $ q_1 ( \tilde{x}_2 ) $
such that
$ p ( \tilde{x}_1 , \tilde{x}_2 ) =   ( \tilde{x}_1 - \theta_{12}t )  
p_1 ( \tilde{x}_1 , \tilde{x}_2 ) +  ( \tilde{x}_2 - \theta_{21}t )
q_1 ( \tilde{x}_2 ) $.
Considering the rest of the points of the second column in the same way as
 what we concluded for  $ q ( \tilde{x}_2 ) $, we see that
$  q_1 ( \tilde{x}_2 ) = ( \tilde{x}_2 - \theta_{22}t )
( \tilde{x}_2 - \theta_{23}t )
 ...( \tilde{x}_2 - \theta_{2m_2}t )  q_1' ( \tilde{x}_2 ) $
for some polynomial $ q_1' ( \tilde{x}_2 )$.
One can now write $ \alpha$ as
$$ \alpha( \tilde{x}_1, \tilde{x}_2) = ( \tilde{x}_1 - \theta_{11}t )   
( \tilde{x}_1 - \theta_{12}t ) p_1 ( \tilde{x}_1, \tilde{x}_2 )
+   ( \tilde{x}_1 - \theta_{11}t )  ( \tilde{x}_2 - \theta_{21}t ) 
... ( \tilde{x}_2 - \theta_{2m_2}t ) q_1' ( \tilde{x}_2 ) $$
$$ + ( \tilde{x}_2 - \theta_{21}t ) ...
( \tilde{x}_2 - \theta_{2m_1}t ) q' ( \tilde{x}_2 ) $$
Proceeding by induction we write $ \alpha $ as
$$ \alpha( \tilde{x}_1, \tilde{x}_2) =  
( \tilde{x}_1 - \theta_{11}t )   ( \tilde{x}_1 - \theta_{12}t )
... ( \tilde{x}_1 - \theta_{1l}t ) q_l $$
$$+ ( \tilde{x}_1 - \theta_{11}t )  ( \tilde{x}_1 - \theta_{12}t ) 
... ( \tilde{x}_1 - \theta_{1(l-1)}t )
( \tilde{x}_2 - \theta_{21}t ) ... 
( \tilde{x}_2 - \theta_{2m_l}t ) q_{l-1}'
+ ... $$
$$ +   ( \tilde{x}_1 - \theta_{11}t )  
( \tilde{x}_2 - \theta_{21}t ) ... 
( \tilde{x}_2 - \theta_{2m_2}t ) q_1'
+ ( \tilde{x}_2 - \theta_{21}t ) ... 
( \tilde{x}_2 - \theta_{2m_1}t ) q'.  $$
This not only completes the proof for the case $ m = 2 $,
but also gives a specific list of the generators in the form 
that was claimed.
  
Inductively assume the lemma is true for any polynomial in
$ m-1 $ variables and for any value of $ c $ so that any two
linear terms in each of the contributing products are distinct. 
We then show it also holds for any polynomial in $ m $ variables
and any value of $ c $ so that any two linear terms in each of 
the contributing products are distinct.

Assume $ \alpha \in \mathcal{I}_+ $ so that it vanishes at the 
given points of $ \CC [t]^m $.
As before we arrange the points at which $ \alpha $ vanishes in 
the following way: the first column consists of the points
$ ( \theta_{11}t, \theta_{2j_2}t ,  ... \theta_{mj_m}t ):
 \quad ( j_2, ..., j_m ) \in A_1 \subset \{ 1,...,n_2 \} 
\times ... \times \{ 1,...,n_m \}, $
the second column is  
$ ( \theta_{12}t, \theta_{2j_2}t ,  ... \theta_{mj_m}t ) $ for
$ ( j_2, ..., j_m ) \in A_2,  $
and the last column is
$ ( \theta_{1l}t, \theta_{2j_2}t ,  ... \theta_{mj_m}t ) $ for
$ ( j_2, ..., j_m ) \in A_l  $
where $ A_i $ are specified by the definition of 
the set $  \bar{\mathcal{F}}_+ $ so that
$ A_i = \{ ( j_2, ..., j_m ) : \quad j_i \leq n_i \textup{ and }
\theta_{1i}  + \sum_{k=2}^m \theta_{kj_k} > c \}$.
Here $l$ is the largest integer such that there is some
point in $ \bar{\mathcal{F}}_+ $ whose first coordinate is 
$ \theta_{1l}t $.

Because 
$ \theta_{11} > \theta_{12} > ... > \theta_{1l} $,
we see that if we have
$\theta_{1(i+1)}  + \sum_{k=2}^m \theta_{kj_k} > c $,
then we also have
$\theta_{1i}  + \sum_{k=2}^m \theta_{kj_k} > c $. Therefore
$ A_l \subset A_{l-1} \subset ... \subset A_2 \subset A_1. $

Considering the first point of the first column we see that 
by the division algorithm
$ \alpha ( \tilde{x}_1, ... , \tilde{x}_m )
= (  \tilde{x}_1 -  \theta_{11}t ) p ( \tilde{x}_1, ... , \tilde{x}_m ) 
+ q  ( \tilde{x}_2, ... , \tilde{x}_m ), $
for some polynomials $p$ and $q$. 
Considering the rest of the points of the first column we find
that $q$ has to satisfy
$ q ( \theta_{2j_2}t ,  ... \theta_{mj_m}t ) = 0 \quad
\textup{ for all } (j_2, ... , j_m ) \in A_1. $
Consider the points in $ \CC [t]^{m-1} $ corresponding to $ A_1 $.
Then 
$ \theta_{2j_2} + ... + \theta_{mj_m} > - \theta_{11} + c  $ 
so that we can apply the induction hypothesis to $ q $ and
$ c' = c - \theta_{11} $ and 
conclude that $ q $ can be written as a linear combination of products
of linear terms of the form
$   ( \tilde{x}_2 -  \theta_{2j_2}t ) $,
$   ( \tilde{x}_3 -  \theta_{3j_3}t ) $, ... and 
$   ( \tilde{x}_2 -  \theta_{mj_m}t ) $,
where $ j_2 $, ..., $ j_m$ are specified by $ A_1 $ and no linear 
term appears twice in each resulting product.

Next considering the second column, we find polynomials
$ p_1 $ and $ q_1 $ such that
$  p ( \tilde{x}_1, ... , \tilde{x}_m ) = 
 (  \tilde{x}_1 -  \theta_{12}t ) p_1  ( \tilde{x}_1, ... , \tilde{x}_m )
+ q_1  ( \tilde{x}_2, ... , \tilde{x}_m ). $
Because $ A_2 \subset A_1 $, we see that
$ q_1 ( \theta_{2j_2}t ,  ... \theta_{mj_m}t ) = 0 \quad
\textup{ for all } (j_2, ... , j_m ) \in A_2. $
So $ q_1 $ is a combination of products of linear terms by the
induction hypothesis so that no linear term appears twice in 
each resulting product.
One can now write $ \alpha $  as
$$ \alpha = ( \tilde{x}_1 - \theta_{11} ) ( \tilde{x}_1 - \theta_{12} )
p_1    +  (\tilde{x}_1 - \theta_{11} ) q_1 + q. $$
Note that the term $ ( \tilde{x}_1 - \theta_{11} ) $
does not appear anywhere in
$ q_1 ( \tilde{x}_2, ... , \tilde{x}_m ) $
so that after multiplying by it in each term of $ q_1 $, 
the linear terms that appear in each resulting product are still
mutually distinct.

Proceeding inductively on the columns we obtain polynomials 
$q$, $q_1$, $q_2$, ... all of which are combinations of products
of linear terms, so that eventually $ \alpha $ can also be written 
in this way with the same property that the linear terms in each
resulting product are mutually distinct.

If the condition in the definition of $ \bar{\mathcal{F}}_+ $ is 
$ \sum_{i=1}^m \theta_{ij_i} < c $, we simply need to write the
table of the points of $ \bar{\mathcal{F}}_- $ in reverse order so 
that the points corresponding to the largest indices appear at
the top of the table. 
Then because
$ \theta_{11} > \theta_{12} > ... > \theta_{1l} $, 
we see that for each $i$ if
$\theta_{1i}  + \sum_{k=2}^m \theta_{kj_k} < c $, then also 
$\theta_{1(i+1)}  + \sum_{k=2}^m \theta_{kj_k} < c $.
Therefore
$ A_1 \subset A_{2} \subset ... \subset A_l. $
We then need to start the argument from the index $ l $ 
proceeding down to $ 1 $.
The rest of the proof is the same.
This finishes the proof of the lemma.
$\hspace{1cm} \stackrel{\infty}{\smile}$

Let us return to geometry and the case of the product of manifolds. 
We shall give a specific representation of some generating sets
of $ \bar{K}_+ $ and $ \bar{K}_- $ which are of the specific
form described in Lemma \ref{algebraiclemma2}.

For simplicity and convenience we relabel the fixed point set
in the following way:

Consider $ \mathcal{A} :=  N_1 \times N_2 \times ... \times N_m, $
where $ N_i =  \{1,...,n_i\} $.
Then we have a one-to-one correspondence between the components
of the fixed point set on which the value of the moment map is the 
same and the elements of $\mathcal{A}$:
$$  \mathbf{F} = \mathbf{F}( J )= ( F_{1j_1}, F_{2j_2}, ... , 
F_{mj_m} ) \sim  J  = ( j_1, j_2, ... , j_m ) \in \mathcal{A}  $$

\begin{defi}
\label{longshortelement}
We define the long elements of  $\mathcal{A}$ as members of the set
$  \mathcal{L} = \{ J \in  \mathcal{A} \mid \mu (\mathbf{F}( J )) > 0 \}, $
and short elements as members of the set
$  \mathcal{S} = \{ J \in  \mathcal{A} \mid \mu (\mathbf{F}( J )) < 0 \}. $
\end{defi}
Consider the projections
$   \left \{\begin{array}{ll}
 P_i : \mathcal{A} \longrightarrow N_i  & \mbox{$ 1 \leq i \leq m$ } \\
 P_i ( j_1, ..., j_i,...,j_m ) = j_i   \end{array}   \right.  $

\begin{defi}
We call a collection
$ \{ \mathcal{A}_i \}_{1 \leq i \leq m} $ where
$  \mathcal{A}_i \subset N_i $
a covering of $ \mathcal{L} \qquad $ \\
(respectively, $ \mathcal{S}$ ), if
$$  \mathcal{L} \subset \cup_{i = 1}^m P_i^{-1} (  \mathcal{A}_i ) \quad
\textup{ (respectively} \quad  \mathcal{S} \subset \cup_{i = 1}^m P_i^{-1}
(  \mathcal{A}_i )). $$
We call it a minimal covering if, whenever we drop just one element
from one of the $ \mathcal{A}_i $'s, it will no longer be a covering
of $\mathcal{L}$ (respectively, $ \mathcal{S}$).

\end{defi}

Notice that some of the $ \mathcal{A}_i $'s may be empty sets and
$ P^{-1}_i( \mathcal{A}_i) = N_1 \times ... \times N_{i-1} \times
\mathcal{A}_i\times N_{i+1} \times ... \times N_m. $

Consider the composition of the map
$$ \CC [t, \tilde{x}_1,...,\tilde{x}_m ]  \stackrel{\eta} \longrightarrow
 H^*_T(M, \CC )  $$
with
$$  H^*_T(M, \CC ) \stackrel{\kappa} \longrightarrow H^*(M_{red}, \CC ). $$
Let $ \bar{K}_+ $ and $ \bar{K}_- $ denote the preimages under $\eta$
of $K_+$ and $K_-$ defined in Theorem \ref{tw1}.
Now we are ready to state our main result:

\begin{theorem}
\label{m}

Consider the case of products of compact connected manifolds such
that the cohomology of each of them is generated by a degree two
form.

(i) The following family of classes of equivariant forms belongs
to and generates $ \bar{K}_+$:
$$  \prod_{1 \leq i \leq m} \prod_{ j_i \in  \mathcal{A}_i }
(\tilde{x}_i - \theta_{ij_i} t ) \hspace{8cm} (1) $$
for all minimal coverings
$\{  \mathcal{A}_i\}$  of $ \mathcal{L} $.

(ii) The following family of classes of equivariant forms belongs
to and generates $ \bar{K}_-$:
$$ \prod_{1 \leq i \leq m}  \prod_{  j_i \in  \mathcal{A}_i }
(\tilde{x}_i - \theta_{ij_i} t ) \hspace{8cm} (2) $$
for all minimal coverings
$ \{ \mathcal{A}_i\} $  of $ \mathcal{S} $.

\end{theorem}

{ \bf Remark:}
The minimality condition was added to avoid some extra terms which do
not contribute to generating $ \bar{K}_+$ or $ \bar{K}_-$.

{ \bf Proof:}
By Lemma \ref{algebraiclemma2}, $ \bar{K}_+ $ has a set of generators
that are products of distinct linear terms.
Moreover the lemma precisely specifies these linear terms :
$ \tilde{x}_i - \theta_{ij_i} t $, where the $  \theta_{ij_i} t$ are the
components of the points in $ \CC [t]^m $ at which the elements of
$  \bar{K}_+ $ vanish.
There are a finite number of polynomials that can be written in this form.
Considering all possible choices there are a total of 
$ 2^{n_1n_2 ... n_m} $ polynomials made out of these linear terms so that 
no linear term appears more than once.
So we have the task of separating all those that belong to $ \bar{K}_+ $
and giving an adequate set of generators for it.

Clearly every element of $ (5) $ vanishes at  $ \mathbf{F} ( J ) $ for all
$ J \in \mathcal{L} $.
Assume $ \alpha \in  \bar{K}_+ $ is a product of the linear terms specified.
We show that $ \alpha $ is a multiple of some polynomial in the class 
$ (5) $. This means the class of polynomials $ (5) $ form a generating set
for $  \bar{K}_+ $.

To show this for each $i$, define the sets $ \mathcal{B}_i $ 
as 
$ \mathcal{B}_i = \{ j_i : \quad ( \tilde{x}_i - \theta_{ij_i} t )
\textup{ divides } \alpha \}. $
Then $ \{ \mathcal{B}_i \}_{i \leq m } $ is a covering of $ \mathcal{L} $,
since $ \alpha \in  \bar{K}_+ $, hence  $ \alpha $ should vanish at
$ \mathbf{F} ( J ) $ for all $ J \in \mathcal{L} $. Then
$ \alpha =  \prod_{1 \leq i \leq m} \prod_{ j_i \in  
\mathcal{B}_i } (\tilde{x}_i - \theta_{ij_i} t ). $
This covering does not have to be a minimal one. 
However it is clear that every covering has a minimal sub-covering, 
i.e., a minimal covering
 $ \{ \mathcal{A}_i \}_{i \leq m } $ such that
$ \mathcal{A}_i \subset \mathcal{B}_i $ for each $ i $.
Then the polynomial in $ (i) $ corresponding to this minimal covering
is a divisor of $ \alpha $ so that the classes $ (5) $ corresponding
to minimal coverings suffice to form a generating set for $  \bar{K}_+ $.
This finishes the proof of $ (i) $.

The proof of $ (ii) $ is similar.
$\hspace{1cm} \stackrel{\infty}{\smile}$

\section{Examples}

{\bf Example 1}
Consider the projective space $ M = \CC P^n $
equipped with a circle action with weights $ m_1,...,m_n $
so that
$ g.[z_0:...:z_n] = [g^{m_0}z_0:...:g^{m_n}z_n] $ for 
$ g \in S^1$  and $ [z_0:...:z_n] \in \CC P^n. $
This action is Hamiltonian with the moment map
$  \mu : \CC P^n \longrightarrow R; \hspace{.5cm}
[z_0:...:z_n] \longmapsto
\frac{\sum_i m_i z_i \bar{z}_i}{\sum_i z_i \bar{z}_i}. $
The fixed points of this action are $ F_i = [ 0:...:1:...:0] $ 
where $1$ is in the $i$-th position.

Kalkman \cite{Ka}  used the localization formula to find the
cohomology ring of the symplectic quotient $ \mu^{-1}(0) / S^1 $.
As we show this is a special case of Theorem \ref{m}:

The cohomology ring of $\CC P^n$ is generated by the degree
two class of the symplectic form $ x \in H^2(\CC P^n) $.  
Define $ \tilde{x} $ as before. Also
$ \mathcal{L} = \{ i \mid \mu (F_i) > 0 \} $ and 
$ \mathcal{S} = \{ i \mid \mu (F_i) < 0 \} $. \\
By Theorem \ref{m}, the polynomials 
$ P = \prod_{ i \in  \mathcal{L}} (\tilde{x} -  \mu (F_i)t) $
and
$ Q = \prod_{ i \in  \mathcal{S}} (\tilde{x} -  \mu (F_i)t) $
(families $(1)$ and $(2)$ of Theorem \ref{m} ) generate 
$ \bar{K}_+$ and $ \bar{K}_-$ respectively. 
They correspond to the minimal coverings $ \{  \mathcal{L} \} $ and
$ \{ \mathcal{S} \} $ of $ \mathcal{L}$ and $ \mathcal{S} $ 
respectively.
This result is the content of Theorem 5.2 in \cite{Ka}.
$\hspace{2cm} \stackrel{\infty}{\smile}$

{\bf Example 2}
Now consider the product of two projective spaces
$ M = \CC P^k \times \CC P^l $ with symplectic forms
$ x_1 $ and $ x_2 $ and a circle acting on both with weights
$ m_0,...,m_k $ and $ n_0,...,n_l $ and moment maps $\mu_1$ and 
$\mu_2$.
Suppose we have ordered the weights so that
$  m_0 >  m_1 > ...  > m_k $ and 
$ n_0 >  n_1 > ...  > n_l $. Then
$ N_1 = \{ 0, ... , k \} $ and 
$ N_2 = \{ 0, ... , l \} $
Consider the diagonal circle action on $M$ and assume $0$ is
a regular value of the moment map $ \mu = \mu_1 + \mu_2 $ 
on $M$ so that
$ m_i + n_j \neq 0 $ for all $ 0 \leq i \leq m_k
\textup{ and } 0 \leq j \leq n_l, $
since this is the value of $ \mu $ on the fixed point with 
$1$ in the $i$-th place in $ \CC P^k $ and in the $j$-th place 
in $ \CC P^l $ and $0$ everywhere else. \\
Note that in the notation of Theorem \ref{m},
$ \theta_{1i} = m_i  \textup{and} \theta_{2j} = n_j.$
Following the explanation in the proof of Lemma 
\ref{algebraiclemma2}
(the case $m=2$ in the notation of that lemma) we see that there 
are integers $q$ and $ l_0 ,..., l_q$ (specified by the weights 
of the actions on  $\CC P^k$ and $\CC P^l$) such that 
$ q \leq k $,  $ l \geq l_0  \geq l_1 \geq .... \geq l_k, $
and for  $ 0 \leq i \leq q$,
$ m_{i} + n_{j} > 0 $ for $  0 \leq j \leq l_i. $

By Theorem \ref{m} and in the notation of that theorem, we obtain
the following classes that generate $ \bar{K}_+$:
$$ (\tilde{x}_1 - m_0 t) (\tilde{x}_1 - m_1 t)... (\tilde{x}_1 - m_q t),$$
$$   (\tilde{x}_1 - m_0 t) (\tilde{x}_1 - m_1 t)... (\tilde{x}_1 - m_{q-1} t)
 (\tilde{x}_2 - n_0 t) ...  (\tilde{x}_2 - n_{l_q} t),$$
$$...$$
$$ (\tilde{x}_1 - m_0 t)
(\tilde{x}_2 - n_0 t) ...  (\tilde{x}_2 - n_{l_1} t), $$
$$(\tilde{x}_2 - n_0 t) ...  (\tilde{x}_2 - n_{l_0} t). $$
Likewise we obtain classes of the above form which generate 
$  \bar{K}_-$
with the only difference that now for $ q' \leq i \leq k $,
$ m_{i} + n_{j} < 0 $ and  for $ l'_i \leq j \leq l, $
for some  $q'$ and $l'_{q'}, ..., l'_l$ that are specified by 
the weights. Then
$ H^*(M_{red}) \cong  \CC [t, \tilde{x}_1, \tilde{x}_2 ]/
\mathcal{I},$
where $ \mathcal{I} $ is the ideal generated by the two 
families of classes introduced  in the example.
$ \quad \stackrel{\infty}{\smile}$

{\bf Example 3} As the next example we would like to consider 
the case of the product of $m$ spheres of radii 
$ r_1, ..., r_m$ and the diagonal circle action.
The result for this case was first obtained by Hausmann and 
Knutson \cite{HK}. They however had a different approach.

So $ M =  S^2_{r_1}\times ... \times  S^2_{r_m} $ and $x_j$ is
the symplectic form of the $j$-th sphere.
The group $G$ is $SU(2)$ or $SO(3)$ acting diagonally on M and 
$T=U(1)$ is its maximal torus acting by rotation around a fixed 
axis, say the $z$-axis on each sphere.
The fixed point set of the circle action on $M$ is then
$ \mathcal{F} = \{(i_1r_1 \hat{k}, ... ,i_mr_m \hat{k}) 
\mid i_j= \pm1 ; 1\leq j \leq m \}, $
where $ \hat{k}$ is the unit vector in the $z$-axis direction.
The moment map of the j-th sphere is
$ \mu_j: S^2_{r_j} \rightarrow \RR  ;
\quad \mu_j(x_j,y_j,z_j)=z_j. $

We label the fixed point set in the following way: Let
$\mathcal{A}= \{ 1,...,m \} $ and $J \subset \mathcal{A}$
an arbitrary subset and consider the following fixed point 
associated to $J$,
$ F_J = (i_1r_1\hat{k}, ... ,i_mr_m\hat{k}),$
where $i_j = 1$ if  $j \in J$ and $i_j = -1 $ if $j \notin J$.

The restriction of each $\tilde{x}_j $ to $F_J$ is given by
\begin{displaymath}
 \tilde{x}_j \mid_{F_J} = \left\{ \begin{array}{ll}
r_j t  & \textrm{ if $j \in J$} \\
-r_jt  & \textrm{if $j \notin J,$}
 \end{array}   \right. 
\end{displaymath}
and the value of the moment map $\mu$ at $F_J$ is
$ \mu (F_J)= \sum_{j \in J}^{} r_j -  \sum_{j \notin J}^{} r_j.$
We assume $0$ is a regular value of the moment map so that 
$ \mu(F_J) \neq 0 $ for all fixed 
points $ F_J $.

\begin{defi}
\label{longshortset}
The set $ J \subset \mathcal{A} $ is called long if $ \mu (F_J) > 0 $,
otherwise it is called short.
The set of all long subsets of $\mathcal{A}$ is denoted by 
$\mathcal{L}$, and that of short subsets is denoted by 
$\mathcal{S}$.
\end{defi}

Therefore $ J \in \mathcal{A} $ is long if and only if
$ \sum_{j \in J}r_j > \sum_{j \notin J}r_j $.

For every subset $J$ of $ \mathcal{A}$ define
$ P_J = \prod_{j \in J} ( \tilde{x}_j - r_j t) $ and
$ Q_J = \prod_{j \in J} ( \tilde{x}_j + r_j t) $
in the equivariant cohomology ring of M.
These are polynomials in the variables $ \tilde{x}_j$.
Consider the following families of classes of polynomials in 
$ H^*_{S^1}(M)$:

$  (i) \hspace{.6cm} ( \tilde{x}_j - r_j t)( \tilde{x}_j + r_j t)
\hspace{3.9cm} j \in \mathcal{A}$

$ (ii)  \hspace{.6cm} P_J \hspace{6.5cm} J \subset \mathcal{A} 
\quad long \hspace{3cm} (3) $

$ (iii)  \hspace{.5cm} Q_J \hspace{6.4cm} J \subset \mathcal{A} 
\quad long $

\begin{theorem}
\label{s}
Let $ M = \prod_i S^2_{r_i} $ and $ \bar{K}_+ $ and $ \bar{K}_- $ be the
preimages under $\eta$ in  $\CC [t, \tilde{x}_1,...,\tilde{x}_m ]$ 
defined above. Then

(a) The families (i) and (ii) together form a set of generators of
$ \bar{K}_+$.

(b)  The families (i) and (iii) together form a set of generators of
$ \bar{K}_-$.

\end{theorem}

\begin{cor}
The cohomology ring of $ M_{red} $ can be written as
$$ H^*(M_{red}) \cong   \CC [t, \tilde{x}_1,...,\tilde{x}_m ]/
\mathcal{I},$$
where $ \mathcal{I} $ is the ideal generated by the families 
(i), (ii) and (iii) in (2).

\end{cor}

{\bf Proof of Theorem \ref{s}:}
Let
$ N_i = \{1,2\} $, $ \theta_{i1} = -r_i $, $ \theta_{i2} = r_i $
and the long/short subsets defined in Definition \ref{longshortset} 
correspond 
to the long/short elements defined in Definition \ref{longshortelement}.
Fix $ 1 \leq j \leq m $ and define
$  \mathcal{A}_j = N_j = \{1,2\} $ and
$  \mathcal{A}_i = \emptyset $ if $ i \neq j $. Then
$ \mathcal{A} = N_1 \times ... \times N_m \subset P_j^{-1}
( \mathcal{A}_j) $, hence
$ \{  \mathcal{A}_i\}_{1 \leq i \leq m} $ is a covering of both
$\mathcal{L}$ and   $\mathcal{S}$, clearly a minimal one in the
notation of Theorem \ref{m}.
The classes $(1)$ and $(2)$ in Theorem \ref{m} corresponding to 
this minimal covering are both
$ ( \tilde{x}_j - r_j t)( \tilde{x}_j + r_j t) $
which is $(i)$ in the collection $(3)$.

Next suppose $L$ is a long element of $\mathcal{A}$, define
$ \mathcal{A}_i = \{2\}$ if $ P_i(L) = 2 $ and $ \emptyset $
otherwise.
To proceed we need to show that any two long subsets have 
nonempty intersection.
In fact if $ J $ and $ L $ are long and 
$ J \cap L = \emptyset $, then
$ \sum_{j \in J}r_j > \sum_{j \notin J}r_j $ and
$ \sum_{j \in L}r_j > \sum_{j \notin L}r_j $ and therefore
$$\sum_{j \in L}r_j > \sum_{j \notin L}r_j =
\sum_{j \in J}r_j + \sum_{j \notin J \cup L}r_j \geq
\sum_{j \in J}r_j > \sum_{j \notin J}r_j =
\sum_{j \in L}r_j + \sum_{j \notin J \cup L}r_j >
\sum_{j \in L}r_j, $$
which is a contradiction. Consequently
$ \{  \mathcal{A}_i\}_{1 \leq i \leq m} $
is a covering of $ \mathcal{L} $.
If $J$ is another long element, there is some $i$ 
such that $ P_i(J) = P_i(L) = \{ 2 \} $, hence
$ J \in P_i^{-1} ( \mathcal{A}_i ) $.
The corresponding class in $(1)$ is then the class $ P_L $
in the collection $(3)$ where here $L$ denotes the long subset 
corresponding to the long element being considered.
If $S$ is a short element, using its long counterpart $L$ 
(in terms of the subsets $ L = \mathcal{A} - S $), we obtain the 
class  $ Q_L $ in $( 3, iii )$.
$\hspace{2cm} \stackrel{\infty}{\smile}$

The elements of the family $ (iii) $ in $(3)$ look different from 
those of the third family introduced in Theorem 6.4 in \cite{HK}. 
They are the same when the coefficient ring 
is $ \CC $. To see this start from the families  $ (3) $ 
in Theorem \ref{s} 
and write $ u_j =  \tilde{x}_j/r_j$. 
The families $(i)$ and $(ii)$ can be written in terms of
$u_j$:
 $$ (\tilde{x}_j + r_j t) (\tilde{x}_j - r_j t)
   = r_j^2 (u_j + t) (u_j - t) $$
and
$$ P_L =  \prod_{j \in L} ( \tilde{x}_j - r_j t)
= (\prod_{j \in L} r_j)  \prod_{j \in L} (\tilde{x}_j/r_j - t)
=   \lambda_L \prod_{j \in L} ( u_j - t) $$
where  $\lambda_L = (\prod_{j \in L}r_j)$.
Every $ Q_L$ for $ L \in \mathcal{L} $ can be rewritten as
$$ Q_L = \prod_{j \in L} ( \tilde{x}_j + r_j t)
 =   \lambda_L \prod_{j \in L} ( u_j + t)
 = \lambda_L \prod_{j \in L} ( u_j - t + 2t ) 
 =  \lambda_L \sum_{J \subset L} \prod_{j \in J} ( u_j - t )
  ( 2t )^{ \mid L-J \mid }. $$
But the long subsets of $L$ have already been included in the 
second 
family $ (ii) $, hence we can drop the terms corresponding to 
$ J \subset L$, $ J \in \mathcal{L} $ in the last expression to 
obtain the classes
$ \lambda_L \sum_{S \subset L, S \in \mathcal{S}}  \quad
\prod_{j \in S} ( u_j -  t ) ( 2t )^{ \mid L-S \mid }. $ 
After dropping the scalar multiples $r_j^2$ and $ \lambda_L $
the new families still generate $ \bar{K}_+$ and $\bar{K}_-$. 

The families introduced in Theorem \ref{s} still do not perfectly 
match with those in Theorem 6.4 in \cite{HK} which in fact are 
finer than ours.

We need to extend our notation: 
let $ \mathcal{A}_m = \{ 1, ...,m \} $. 
Consider $r=r_m$ in Theorem \ref{tw1}. Define 
$ \mathcal{L}(r_m) = \{ L \subset  \mathcal{A}_{m-1} : 
\sum_L r_j - \sum_{\mathcal{A}_{m-1} - L} r_j > r_m \}, $ 
and
$ \mathcal{S}(r_m) = \{ S \subset  \mathcal{A}_{m-1} : 
\sum_S r_j - \sum_{\mathcal{A}_{m-1} - S} r_j < r_m \}. $
Consider the following families: 
$$ \hspace{-2.6cm}(i)' \hspace{.9cm} (u_j -  t)
( u_j +  t) \hspace{6.8cm} j \in \mathcal{A}_{m-1} $$
$$  \hspace{-2.1cm}(ii)' \hspace{.8cm} P'_J =
 \prod_{j \in L} ( u_j - t) \hspace{7cm} 
J  \in \mathcal{L}(r_m) $$
$$  \hspace{-2cm}(iii)' \hspace{.7cm} Q'_J
= \sum_{S \subset L, S \in \mathcal{S}}  \quad \prod_{j \in S} 
( u_j -  t ) ( 2t )^{ \mid L-S \mid } \hspace{3.8cm}
J \in \mathcal{L}(r_m), $$
where  $ u_j =  \tilde{x}_j/r_j$ and
$ \lambda_L = \prod_{j \in L}r_j $ as in remark (1).
Then $r_m$ is a regular value of the moment map for the abelian
polygon space which is defined as 
$$ M = \prod_{i=1}^{m-1} S^2_{r_i} //_{r_m} SO_2. $$

\begin{cor}
\label{s2}
Consider $ \bar{K}_+(r_m)$ and $ \bar{K}_-(r_m)$ for the abelian 
polygon space. Then \\
(a) The families $(i)'$ and $(ii)'$ together form a set of generators 
of $ \bar{K}_+(r_m)$.\\
(b)  The families $(i)'$ and $(iii)'$ together form a set of generators 
of $ \bar{K}_-(r_m)$.

\end{cor}

The proof is the same as that of Theorem \ref{s} adding the 
comments that we gave after the proof of that theorem.
$\hspace{1cm} \stackrel{\infty}{\smile}$

The following corollary is Theorem 6.4 in \cite{HK} when the 
coefficient ring is $ \CC $:

\begin{cor} 
\label{hk6.4}
 For the abelian polygon space,  
$ \bar{K}(r_m)$ is generated by the families $(i)''$, $(ii)''$ and 
$(iii)''$
 which are defined as follows: 
$$  \hspace{-2.4cm} (i)'' \hspace{.7cm} ( u_j -  t)
( u_j +  t)
\hspace{6cm} j \in \mathcal{A}_{m-1} $$
$$  \hspace{-2.7cm} (ii)'' \hspace{.5cm} P''_L = \prod_{j \in L} 
( u_j - t) \hspace{5.8cm}  L \in \mathcal{L}_m $$
$$ \hspace{-1cm} (iii)''  \hspace{.5cm} Q''_L = \sum_{S \subset L, 
S \in \mathcal{S}_m} \prod_{j \in S} ( u_j - t) 
 (2t)^{\mid L - S \mid}  \hspace{2.8cm} L  
\in \mathcal{P}(\mathcal{A}_{m-1}) \cap  \mathcal{L}.$$
Here $\mathcal{P}(\mathcal{A}_{m-1}) $ is the set of all subsets 
of $\mathcal{A}_{m-1}$, and
 $ \mathcal{L}_m $ and  $ \mathcal{S}_m $ are defined as 
$  \mathcal{L}_m = \{ L \subset \mathcal{A}_{m-1} : \quad L \cup 
\{m\} \in \mathcal{L} \}, $ and 
$ \mathcal{S}_m = \{ S \subset \mathcal{A}_{m-1} : \quad S 
\cup \{m\} \in \mathcal{S} \}. $

\end{cor}

{\bf Proof:}
The argument used after the proof of Theorem \ref{s} 
applies here too to show that the 
classes $Q_J$ in $(iii)'$ can be replaced by
$ (iii)''' \hspace{.2cm}
\sum_{S \subset L, S \in \mathcal{S}(r_m)}  
\prod_{j \in S} ( u_j -  t )
( 2t )^{ \mid L-S \mid }  $ for 
 $ J \in \mathcal{L}(r_m) $,
which together with $(i)'$ and $(ii)'$ still generate 
$ \bar{K}(r_m)$.

Then, notice that
$ \mathcal{L}(r_m) \subset  \mathcal{L} \cap \mathcal{P}
(\mathcal{A}_{m-1}) \subset \mathcal{L}_m $,
hence the families $ (ii)''$ and $ (iii)''$ are larger than the
families  $ (ii)'$ and $ (iii)'''$ respectively.
Furthermore this allows us to remove some of the terms in the elements
in $ (iii)'''$ to obtain the elements in  $ (iii)''$. 
In fact consider a term in a class in  $ (iii)'''$ corresponding to some
$ S \in \mathcal{S}(r_m) - \mathcal{S}_m $. 
This means
$ S \cup \{m \} \in \mathcal{L} $, hence $ S \in \mathcal{L}_m $. 
Thus the term corresponding to this $S$ has already been considered
in  $ (ii)''$. 
This establishes that the new families suffice to generate  
$ \bar{K}(r_m). \hspace{1cm} \stackrel{\infty}{\smile}$

{\bf Remark:}
 The classes $ V_j$ and $R$ in Theorem 6.4 in \cite{HK} correspond 
to our classes $ u_j - t $ and $ 2t$ respectively.

\section{ Acknowledgment }
First of all, I would like to thank my advisor Professor Lisa 
Jeffrey for her
great support and patience. She has been of great help to me.

Then I would like to thank professors R. Buchweitz, R. Goldin, 
E. Meinrenken, M. Spivakovsky, and S. Tolman for the very helpful 
discussions that I had with them.
I also thank  Professor G. Elliott and my friends Kiumars Kaveh 
and Lila Rasekh that have been supportive to me during my 
graduate studies.

\vspace{1cm}

Ramin Mohammadalikhani \\
Department of Mathematics \\
University of Toronto \\
Toronto, Ontario \\
Canada, M5S-3G3 \\


\thebibliography{bibliography}

\bibitem [Au]{Au} M. Audin
{\it The topology of torus actions on symplectic manifolds},
Birkhauser (Progress in Mathematics v.{ \bf 93}) 1991.

\bibitem [BGV]{BGV} N. Berline, E. Getzler, M. Vergne 
{\it Heat Kernels and Dirac Operators},
Springer-Verlag, 1992.

\bibitem [G]{G} R. Goldin {\it The cohomology of weight 
varieties}, Ph.D. thesis, MIT, 1999.

\bibitem [G2]{G2} R. Goldin
``The cohomology of weight varieties and polygon spaces''
{\it Advances in Mathematics} {\bf 160} (2001)  175-204.

\bibitem [HK]{HK} J. Hausmann, A. Knutson
``The cohomology ring of polygon spaces''
{\it Ann. Inst. Fourier} {\bf 48} (1998) 281-321.

\bibitem [JK1]{JK1} L. Jeffrey, F. Kirwan
 ``Localization for nonabelian group actions''
{ \it Topology} { \bf 34} (1995) 291-327.

\bibitem [JK2]{JK2} L. Jeffrey, F. Kirwan
``Intersection theory on moduli spaces of
holomorphic bundles of arbitrary rank on a Riemann surface''
{ \it Annals of Math. } {\bf 148 } (1998) 109-191.

\bibitem [JK3]{JK3}  L. Jeffrey, F. Kirwan 
``Localization and the quantization conjecture''
{ \it Topology} {\bf 36} (1997) 647-693.

\bibitem [Ka]{Ka} J. Kalkman 
``Cohomology rings of symplectic quotients''
{\it J. Reine Angew. Math.} {\bf 458} (1995) 37-52.

\bibitem [Ki1]{Ki1} F. Kirwan
{\it Cohomology of Quotients in Symplectic and Algebraic Geometry},
Princeton University Press, 1984. 

\bibitem [TW1]{TW1} S. Tolman, J. Weitsman
``The cohomology rings of symplectic quotients''
Preprint; {\it Commun. in Analysis and Geometry}, to appear.

\bibitem [TW2]{TW2}  S. Tolman, J. Weitsman
``On semifree symplectic circle actions with isolated fixed points''
{\it Topology} {\bf 39} (2000) 299-309.
\end{document}